\documentclass[twoside]{article}
\usepackage{amsmath,amssymb,amscd,amsthm,verbatim,alltt,amsfonts,array}
\usepackage{mathrsfs}
\usepackage[english]{babel}
\usepackage{latexsym}
\usepackage{amssymb}
\usepackage{euscript}
\usepackage{graphicx}

 \def\NZQ{\Bbb}               
 \def\NN{{\NZQ N}}
 
 \def\ZZ{{\NZQ Z}}

 \def\C{{\mathcal C}}

\def\Xb{{\bold X}}
 \def\ab{{\bold a}}
 
 \def\xb{{\bold x}}

 \def\opn#1#2{\def#1{\operatorname{#2}}} 
 \opn\chara{char} \opn\length{\ell} \opn\pd{pd} \opn\rk{rk}
 \opn\projdim{proj\,dim} \opn\injdim{inj\,dim} \opn\rank{rank}
 \opn\depth{depth} \opn\grade{grade} \opn\height{height}
 \opn\embdim{emb\,dim} \opn\codim{codim}
 
 \opn\Tr{Tr} \opn\bigrank{big\,rank}
 \opn\superheight{superheight}\opn\lcm{lcm}
 \opn\trdeg{tr\,deg}
 \opn\reg{reg} \opn\lreg{lreg} \opn\ini{in} \opn\lpd{lpd}
 \opn\size{size} \opn\sdepth{sdepth}
 \opn\link{link}\opn\fdepth{fdepth}\opn\lex{lex}
 \opn\GCD{GCD}
 \opn\div{div} \opn\Div{Div} \opn\cl{cl} \opn\Cl{Cl}
 \opn\Spec{Spec} \opn\Supp{Supp} \opn\supp{supp} \opn\Sing{Sing}
 \opn\Ass{Ass} \opn\Min{Min}\opn\Mon{Mon}
 \opn\astab{astab}
 \opn\Ann{Ann} \opn\Rad{Rad} \opn\Soc{Soc}
 \opn\Im{Im} \opn\Ker{Ker} \opn\Coker{Coker} \opn\Am{Am}
 \opn\Hom{Hom} \opn\Tor{Tor} \opn\Ext{Ext} \opn\End{End}
 \opn\Aut{Aut} \opn\id{id}
 
 \opn\nat{nat}
 \opn\pff{pf}
 \opn\Pf{Pf} \opn\GL{GL} \opn\SL{SL} \opn\mod{mod} \opn\ord{ord}
 \opn\Gin{Gin} \opn\Hilb{Hilb}\opn\sort{sort}
 \opn\Tot{Tot}
 \opn\astab{astab}
 \opn\aff{aff} \opn
 \con{conv} \opn\relint{relint} \opn\st{st}
 \opn\lk{lk} \opn\cn{cn} \opn\core{core} \opn\vol{vol}
 \opn\link{link} \opn\star{star}\opn\lex{lex}\opn\set{set}
 \opn\gr{gr}
 
 \def\pot#1#2{#1[\kern-0.28ex[#2]\kern-0.28ex]}

 \opn\dirlim{\underrightarrow{\lim}}
 \opn\inivlim{\underleftarrow{\lim}}
 \let\union=\cup

 \def\Implies{\ifmmode\Longrightarrow \else
         \unskip${}\Longrightarrow{}$\ignorespaces\fi}
 \def\implies{\ifmmode\Rightarrow \else
         \unskip${}\Rightarrow{}$\ignorespaces\fi}
 \def\iff{\ifmmode\Longleftrightarrow \else
         \unskip${}\Longleftrightarrow{}$\ignorespaces\fi}

 \let\:=\colon
 \newtheorem{Theorem}{Theorem}[section]
 \newtheorem{Lemma}[Theorem]{Lemma}
 \newtheorem{Corollary}[Theorem]{Corollary}
 \newtheorem{Proposition}[Theorem]{Proposition}
 \newtheorem{Remark}[Theorem]{Remark}
 
 \newtheorem{Example}[Theorem]{Example}
 
 \newtheorem{Definition}[Theorem]{Definition}

 \let\epsilon\varepsilon
 \let\kappa=\varkappa
 \def\qed{\ifhmode\textqed\fi
       \ifmmode\ifinner\quad\qedsymbol\else\dispqed\fi\fi}
 \def\textqed{\unskip\nobreak\penalty50
        \hskip2em\hbox{}\nobreak\hfil\qedsymbol
        \parfillskip=0pt \finalhyphendemerits=0}
 \def\dispqed{\rlap{\qquad\qedsymbol}}
 
 \opn\dis{dis}
 \def\pnt{{\raise0.5mm\hbox{\large\bf.}}}
 
 \opn\Lex{Lex}

\newcommand{\salt}{\vspace*{2.5mm}}     

\def\boxit#1{\vbox{\hrule\hbox{\vrule\kern2pt
 \vbox{\kern4pt#1\kern2pt}\kern2pt\vrule}\hrule}}
\def\qed{\hfill
 $\hskip0.3cm
 \boxit{\hsize 2pt \vsize 10pt}$\bigskip\noindent}

\def\C{{\mathchoice {\setbox0=\hbox{$\displaystyle\rm C$}\hbox{\hbox
to0pt{\kern0.4\wd0\vrule height0.9\ht0\hss}\box0}}
{\setbox0=\hbox{$\textstyle\rm C$}\hbox{\hbox
to0pt{\kern0.4\wd0\vrule height0.9\ht0\hss}\box0}}
{\setbox0=\hbox{$\scriptstyle\rm C$}\hbox{\hbox
to0pt{\kern0.4\wd0\vrule height0.9\ht0\hss}\box0}}
{\setbox0=\hbox{$\scriptscriptstyle\rm C$}\hbox{\hbox
to0pt{\kern0.4\wd0\vrule height0.9\ht0\hss}\box0}}}}
\def\Q{{\mathchoice {\setbox0=\hbox{$\displaystyle\rm
Q$}\hbox{\raise
0.15\ht0\hbox to0pt{\kern0.4\wd0\vrule height0.8\ht0\hss}\box0}}
{\setbox0=\hbox{$\textstyle\rm Q$}\hbox{\raise
0.15\ht0\hbox to0pt{\kern0.4\wd0\vrule height0.8\ht0\hss}\box0}}
{\setbox0=\hbox{$\scriptstyle\rm Q$}\hbox{\raise
0.15\ht0\hbox to0pt{\kern0.4\wd0\vrule height0.7\ht0\hss}\box0}}
{\setbox0=\hbox{$\scriptscriptstyle\rm Q$}\hbox{\raise
0.15\ht0\hbox to0pt{\kern0.4\wd0\vrule height0.7\ht0\hss}\box0}}}}
\def\T{{\mathchoice {\setbox0=\hbox{$\displaystyle\rm
T$}\hbox{\hbox to0pt{\kern0.3\wd0\vrule height0.9\ht0\hss}\box0}}
{\setbox0=\hbox{$\textstyle\rm T$}\hbox{\hbox
to0pt{\kern0.3\wd0\vrule height0.9\ht0\hss}\box0}}
{\setbox0=\hbox{$\scriptstyle\rm T$}\hbox{\hbox
to0pt{\kern0.3\wd0\vrule height0.9\ht0\hss}\box0}}
{\setbox0=\hbox{$\scriptscriptstyle\rm T$}\hbox{\hbox
to0pt{\kern0.3\wd0\vrule height0.9\ht0\hss}\box0}}}}
\def\Z{{\mathchoice {\hbox{$\sf\textstyle Z\kern-0.4em Z$}}
{\hbox{$\sf\textstyle Z\kern-0.4em Z$}}
{\hbox{$\sf\scriptstyle Z\kern-0.3em Z$}}
{\hbox{$\sf\scriptscriptstyle Z\kern-0.2em Z$}}}}

\newcommand{\eqnoone}  
   {}
\newcommand{\eqnotwo}  
   {}

\newcounter{alf}

\newcommand{\adresa}[1]{\par\vspace*{-11pt}
                        \begin{flushright}
                        {\small
                        #1}
                        \end{flushright}
                        }

\begin{document}

  \vskip 1.2 true cm
\setcounter{page}{1}

 \begin{center} {\bf Combinatorial properties of strong quasi-$n$-partite graphs } \\
          \medskip

 { Monica La Barbiera and Roya Moghimipor}
\end{center}

\begin{abstract}
Monomial ideals corresponding to strong quasi-$n$-partite graphs are considered.
Some algebraic and combinatorial properties of generalized graph ideals of a strong quasi-$n$-partite graph are studied.
Furthermore, we show that the edge ideal of a strong quasi-$n$-partite graph is Cohen-Macaulay.
\end{abstract}

\begin{quotation}

\end{quotation}

\thispagestyle{empty}

\section{Introduction}
\label{0ne}
Let $G$ be a graph on vertex set $V$.
Let $K$ be a field and $S=K[x_1,\dots,x_n]$ the polynomial
ring in $n$ variables over $K$ with each $x_i$ of degree $1$.
Next we consider the polynomial ring $T$ over $K$ in the variables
\[
x_{11},\ldots,x_{1m_1},x_{21},\ldots,x_{2m_2},\ldots,x_{n1},\ldots,x_{nm_n}.
\]
In the present paper we consider classes of monomial ideals of $T$ with linear resolution that can arise from graph theory \cite{BM, H}. Algebraic objects attached to $G$ are the edge ideals $I(G)$ and the generalized graph ideals $I_{t}(G)$. The main objective of this paper is to study combinatorial and algebraic properties of generalized graph ideals of a strong quasi-$n$-partite graph.

A graph $G$ is said to be $n$-partite graph if its vertex set $V=V_1\union V_2\union \cdots\union V_n$ and $V_i=\{x_{i1},\ldots,x_{im_i}\}$ for $i=1,\ldots,n$, every edge joins a vertex of $V_{i}$ with a vertex of $V_{i+1}$.
When $n=2$ these are the bipartite graphs.

A quasi-$n$-partite graph is an $n$-partite graph having loops, a strong quasi-$n$-partite graph is a complete $n$-partite graph having loops in all its vertices. A strong quasi-n-partite graph on vertices $x_{11},\ldots,x_{1m_1},\ldots,x_{n1},\ldots,x_{nm_n}$ is denoted by $\mathcal{K}'_{m_1,\dots,m_n}$. In several papers on graph theory some algebraic properties of strong quasi-bipartite graphs are studied (\cite{IB, BPR}).

Matroid theory is a very active and fascinating research area in combinatorics.
The discrete polymatroid is a multiset analogue of the matroid. We refer to \cite{HHC, HH, VI} for the theory of discrete polymatroids.

The class of polymatroidal ideals is one of the rare classes of monomial ideals with the property that all powers of an ideal in this class have a
linear resolution. In fact, the powers of a polymatroidal ideal are polymatroidal ideals and the polymatroidal
ideals have linear quotients; hence, they are an important class of monomial ideals with linear resolution \cite{HT1}.

In \cite{LM}, the authors introduced the generalized notion of a generalized discrete bi-polymatroid and the corresponding notion of a generalized
bi-polymatroidal ideal in $T$.
Let $L^{*}$ be a monomial ideal of $T$ generated in a single degree, and let $G(L^{*})$ be
its unique set of minimal monomial generators. Then $L^{*}$ is called a generalized bi-polymatroidal ideal if for any two monomials
$u = x_{11}^{a_{11}} \cdots x_{1m_1}^{a_{1m_{1}}} \cdots x_{n1}^{a_{n1}} \cdots x_{nm_n}^{a_{nm_{n}}}$
and
$v = x_{11}^{b_{11}} \cdots x_{1m_1}^{b_{1m_{1}}} \cdots x_{n1}^{b_{n1}} \cdots x_{nm_n}^{b_{nm_{n}}}$
in $G(L^{*})$ and for each $ij$ with $a_{ij} > b_{ij}$, then there exists $l\in \{1,\ldots,m_{i}\}$ with $a_{il} < b_{il}$
such that $x_{il}(u/x_{ij}) \in G(L^{*})$.

The present paper is organized as follows. In Section \ref{one} we study monomial ideals arising from strong quasi-$n$-partite graphs.
In Theorem \ref{Cm} we show that the edge ideal of a strong quasi-$n$-partite graph is Cohen-Macaulay.
Furthermore, we compute powers of $I_{t}(\mathcal{K}'_{m_1,\dots,m_n})$, see Proposition \ref{power}.
In Theorem \ref{power-bi-polymatroid} we show that powers of the generalized graph ideal $I_{t}(\mathcal{K}'_{m_1,\dots,m_n})$ associated to a strong quasi-$n$-partite graph $K'_{m_1,\dots,m_n}$ is a generalized bi-polymatroidal ideal.

In Section \ref{two} we study a certain permanence property of $I_{t}(\mathcal{K}'_{m_1,\dots,m_n})$.
In Theorem \ref{colon-generalized} we show that $I_{t}(\mathcal{K}'_{m_1,\dots,m_n}):u$ is a generalized bi-polymatroidal ideal for all monomials $u$ of $T$.

Let $\wp$ be a prime ideal of $T$. The monomial localization of a monomial ideal $I_{t}(\mathcal{K}'_{m_1,\dots,m_n})$ with
respect to $\wp$ is the monomial ideal $I_{t}(\mathcal{K}'_{m_1,\dots,m_n})(\wp)$ obtained from $I_{t}(\mathcal{K}'_{m_1,\dots,m_n})$ by substituting 1 for the variables that do not belong to $\wp$. In Theorem \ref{localization} we prove that monomial localizations of $I_{t}(\mathcal{K}'_{m_1,\dots,m_n})$ at monomial prime ideals is a generalized bi-polymatroidal ideal.

\section{Generalized graph ideals}
\label{one}
Let $K$ be a field and $S=K[x_1,\dots,x_n]$ the polynomial ring in $n$ variables over $K$ with each $x_i$ of degree $1$.
Let $I\subset S$ be a monomial ideal with $I\neq S$ whose minimal set of generators is $G(I)=\{\xb^{\ab_1},\ldots, \xb^{\ab_m}\}$. Here $\xb^{\ab}=x_1^{\ab(1)}x_2^{\ab(2)}\cdots x_n^{\ab(n)}$ for $\ab=(\ab(1),\ldots,\ab(n))\in\NN^n$.

Next we consider the polynomial ring $T$ over $K$ in the variables
\[
x_{11},\ldots,x_{1m_1},x_{21},\ldots,x_{2m_2},\ldots,x_{n1},\ldots,x_{nm_n}.
\]
We want to study  monomial ideals arising from strong quasi-n-partite graphs.

A graph $G$ consists of a finite set $V=\{x_1,\dots,x_n\}$ of vertices and a collection $E(G)$ of subsets of $V$, that consists of pairs $\{x_i,x_j\}$, for some
$x_i,x_j\in V$.
A graph $G$ has loops if it is not requiring $x_i\neq x_j$ for all edges $\{x_i,x_j\}$ of $G$. Thus the edge $\{x_i,x_i\}$ is said a loop of $G$.

\begin{Definition}
\label{quasi-n-partite}
A graph $G$ with loops is said to be {\em quasi-n-partite} if its vertex set $V=V_1\union V_2\union \cdots\union V_n$ and $V_i=\{x_{i1},\ldots,x_{im_i}\}$ for $i=1,\ldots,n$, every edge joins a vertex of $V_{i}$ with a vertex of $V_{i+1}$, and some vertices in $V$ have loops.
\end{Definition}

\begin{Definition}
\label{strong}
A quasi-n-partite graph $G$ is called {\em strong} if it is a complete n-partite graph and all its vertices have loops.
\end{Definition}

A strong quasi-n-partite graph on vertices $x_{11},\ldots,x_{1m_1},\ldots,x_{n1},\ldots,x_{nm_n}$ will be denoted by $\mathcal{K}'_{m_1,\dots,m_n}$.

\begin{Definition}
\label{walk}
Let $G$ be a graph with loops in each of its $n$ vertices. A {\em walk} of {\em length} $t$ in $G$ is an alternating sequence
\[
w=\{v_{i_{0}},l_{i_{1}},v_{i_{1}},l_{i_{2}},\dots, v_{i_{t-1}}, l_{i_{t}},v_{i_{t}}\},
\]
where $v_{i_{0}}$ or $v_{i_{g}}$ is a vertex of $G$ and
$l_{i_{g}}=\{v_{i_{g-1}},v_{i_{g}}\}$, $g=1,\dots,t$, is either the edge joining $v_{i_{g-1}}$ and $v_{i_{g}}$ or a loop if $v_{i_{g-1}}=v_{i_{g}}, 1\leq i_{0}\leq i_{1}\leq \dots\leq i_{t}\leq n$.
\end{Definition}

\begin{Example}
{\em
Let $G$ be a strong quasi-bipartite graph on vertices $x_{1},\ldots,x_{n}$, $y_{1},\ldots,y_{m}$.
A walk of length 2 in $G$ is
\[
\{x_i,l_{i},x_i,l_{ij},y_{j}\} \quad  \text{or} \quad \{x_i,l_{ij},y_j,l_j,y_j\}
\]
where $l_{i}=\{x_i,x_i\}$, $l_{j}=\{y_j,y_j\}$ are loops, and $l_{ij}$
is the edge joining $x_i$ and $y_j$. Because $G$ is bipartite, any walk in it have not the edges $\{x_{i_{h}},x_{i_{k}}\}$, $i_h\neq i_k$, and
$\{y_{j_{h}},y_{j_{k}}\}$, $j_h\neq j_k$.
}
\end{Example}

Let $G$ be a strong quasi-n-partite graph with vertex set $V=V_1\union V_2\union \cdots\union V_n$ and $V_i=\{x_{i1},\ldots,x_{im_i}\}$ for $i=1,\ldots,n$.
The {\em generalized graph ideal} $I_t(G)$ associated to $G$ is the ideal of the polynomial ring $T=K[x_{11},\ldots,x_{1m_1},\dots, x_{n1},\ldots,x_{nm_n}]$ generated by all the monomials of degree $t$ corresponding to the walks of length $t-1$. By definition the ideal $I_{t}(G)$ is not trivial for $2\leq t\leq 2(m_{1}+\dots+m_{n})-1$.

Now let $\mathcal{K}'_{m_1,\dots,m_n}$ be a strong quasi-n-partite graph with vertex set $V=V_1\union V_2\union \cdots\union V_n$ and $V_i=\{x_{i1},\ldots,x_{im_i}\}$ for $i=1,\ldots,n$. Let $t, q_{1}, \dots ,q_{n}$ be non negative integers and $\sum_{i=1}^{n}q_{i}=t$, $q_{1}, \dots ,q_{n}\geq0$.
For this graph we have
\[
I_{t}(\mathcal{K}'_{m_1,\dots,m_n})=\sum_{\sum_{i=1}^{n}q_{i}=t, q_{i}\neq t}L_{1,q_{1},2}\cdots L_{n,q_{n},2},
\]
where the ideals $L_{i,q_{i},2}$ are Veronese type ideals of degree $q_{i}$ generated by the monomials $x_{i1}^{a_{i1}} \cdots x_{im_{i}}^{a_{im_{i}}}$ with $\sum_{j=1}^{m_{i}}a_{ij}=q_{i}$ and $0 \leq a_{ij}\leq 2$ for $i=1,\dots,n$. Here we set $L_{i,q_{i},2}=(0)$, if $q_{i}=t$.
\begin{Remark}
{\em
If $t=2$, the ideal $I_{t}(\mathcal{K}'_{m_1,\dots,m_n})$ does not describe the edge ideal
\[
I(\mathcal{K}'_{m_1,\dots,m_n})=I_{2}(\mathcal{K}'_{m_1,\dots,m_n})
 \]
of a strong quasi-n-partite graph. Let $\mathcal{K}'_{2,2}$ be the strong quasi-bipartite graph on vertices $x_{11},x_{12},x_{21},x_{22}$,
then
\[
I(\mathcal{K}'_{2,2})=(x_{11}x_{21},x_{11}x_{22},x_{12}x_{21},x_{12}x_{22},x_{11}^{2},x_{12}^{2},x_{21}^{2},x_{22}^{2}),
\]
but
$I_{2}(\mathcal{K}'_{2,2})=(x_{11}x_{21},x_{11}x_{22},x_{12}x_{21},x_{12}x_{22})$. Therefore, $I(\mathcal{K}'_{2,2})\neq I_{2}(\mathcal{K}'_{2,2})$.
}
\end{Remark}

\begin{Example}
\label{L9}
{
Let $T=K[x_{11},x_{12},x_{21},x_{22},x_{31},x_{32}]$ be a polynomial ring and $\mathcal{K}'_{2,2,2}$ be the strong quasi-3-partite graph on vertices $x_{11},x_{12},x_{21},x_{22},x_{31},x_{32}$. Then
\begin{eqnarray*}
&I_{9}(\mathcal{K}'_{2,2,2})=&
L_{1,1,2}L_{2,4,2}L_{3,4,2}+L_{1,2,2}L_{2,3,2}L_{3,4,2}+L_{1,2,2}L_{2,4,2}L_{3,3,2}
\\
&&
+L_{1,3,2}L_{2,2,2}L_{3,4,2}+L_{1,3,2}L_{2,3,2}L_{3,3,2}+L_{1,3,2}L_{2,4,2}L_{3,2,2}
\\
&&
+L_{1,4,2}L_{2,1,2}L_{3,4,2}+L_{1,4,2}L_{2,2,2}L_{3,3,2}+L_{1,4,2}L_{2,3,2}L_{3,2,2}
\\
&&
+L_{1,4,2}L_{2,4,2}L_{3,1,2},
\end{eqnarray*}
where the ideals $L_{i,q_{i},2}$ are Veronese type ideals of degree $q_{i}$ generated by the monomials $x_{i1}^{a_{i1}} \cdots x_{im_{i}}^{a_{im_{i}}}$
with $\sum_{j=1}^{m_{i}}a_{ij}=q_{i}$ and $0 \leq a_{ij}\leq 2$ for $i=1,2,3$.
Therefore,
\begin{eqnarray*}
&I_{9}(\mathcal{K}'_{2,2,2})=&
(x_{11}x_{21}^{2}x_{22}^{2}x_{31}^{2}x_{32}^{2},
x_{12}x_{21}^{2}x_{22}^{2}x_{31}^{2}x_{32}^{2},
x_{11}^{2}x_{21}^{2}x_{22}x_{31}^{2}x_{32}^{2},
\\
&&
x_{11}^{2}x_{21}x_{22}^{2}x_{31}^{2}x_{32}^{2},
x_{11}x_{12}x_{21}^{2}x_{22}x_{31}^{2}x_{32}^{2},
x_{11}x_{12}x_{21}x_{22}^{2}x_{31}^{2}x_{32}^{2},
\\
&&
x_{12}^{2}x_{21}^{2}x_{22}x_{31}^{2}x_{32}^{2},
x_{12}^{2}x_{21}x_{22}^{2}x_{31}^{2}x_{32}^{2},
x_{11}^{2}x_{21}^{2}x_{22}^{2}x_{31}^{2}x_{32},
\\
&&
x_{11}^{2}x_{21}^{2}x_{22}^{2}x_{31}x_{32}^{2},
x_{11}x_{12}x_{21}^{2}x_{22}^{2}x_{31}^{2}x_{32},
x_{11}x_{12}x_{21}^{2}x_{22}^{2}x_{31}x_{32}^{2},
\\
&&
x_{12}^{2}x_{21}^{2}x_{22}^{2}x_{31}^{2}x_{32},
x_{12}^{2}x_{21}^{2}x_{22}^{2}x_{31}x_{32}^{2},
x_{11}^{2}x_{12}x_{21}^{2}x_{31}^{2}x_{32}^{2},
\\
&&
x_{11}^{2}x_{12}x_{21}x_{22}x_{31}^{2}x_{32}^{2},
x_{11}^{2}x_{12}x_{22}^{2}x_{31}^{2}x_{32}^{2},
x_{11}x_{12}^{2}x_{21}^{2}x_{31}^{2}x_{32}^{2},
\\
&&
x_{11}x_{12}^{2}x_{21}x_{22}x_{31}^{2}x_{32}^{2},
x_{11}x_{12}^{2}x_{22}^{2}x_{31}^{2}x_{32}^{2},
x_{11}^{2}x_{12}x_{21}^{2}x_{22}x_{31}^{2}x_{32},
\\
&&
x_{11}^{2}x_{12}x_{21}^{2}x_{22}x_{31}x_{32}^{2},
x_{11}^{2}x_{12}x_{21}x_{22}^{2}x_{31}^{2}x_{32},
x_{11}^{2}x_{12}x_{21}x_{22}^{2}x_{31}x_{32}^{2},
\\
&&
x_{11}x_{12}^{2}x_{21}^{2}x_{22}x_{31}^{2}x_{32},
x_{11}x_{12}^{2}x_{21}^{2}x_{22}x_{31}x_{32}^{2},
x_{11}x_{12}^{2}x_{21}x_{22}^{2}x_{31}^{2}x_{32},
\\
&&
x_{11}x_{12}^{2}x_{21}x_{22}^{2}x_{31}x_{32}^{2},
x_{11}^{2}x_{12}x_{21}^{2}x_{22}^{2}x_{31}^{2},
x_{11}^{2}x_{12}x_{21}^{2}x_{22}^{2}x_{31}x_{32},
\\
&&
x_{11}^{2}x_{12}x_{21}^{2}x_{22}^{2}x_{32}^{2},
x_{11}x_{12}^{2}x_{21}^{2}x_{22}^{2}x_{31}^{2},
x_{11}x_{12}^{2}x_{21}^{2}x_{22}^{2}x_{31}x_{32},
\\
&&
x_{11}x_{12}^{2}x_{21}^{2}x_{22}^{2}x_{32}^{2},
x_{11}^{2}x_{12}^{2}x_{21}x_{31}^{2}x_{32}^{2},
x_{11}^{2}x_{12}^{2}x_{22}x_{31}^{2}x_{32}^{2},
\\
&&
x_{11}^{2}x_{12}^{2}x_{21}^{2}x_{31}^{2}x_{32},
x_{11}^{2}x_{12}^{2}x_{21}^{2}x_{31}x_{32}^{2},
x_{11}^{2}x_{12}^{2}x_{21}x_{22}x_{31}^{2}x_{32},
\\
&&
x_{11}^{2}x_{12}^{2}x_{21}x_{22}x_{31}x_{32}^{2},
x_{11}^{2}x_{12}^{2}x_{22}^{2}x_{31}^{2}x_{32},
x_{11}^{2}x_{12}^{2}x_{22}^{2}x_{31}x_{32}^{2},
\\
&&
x_{11}^{2}x_{12}^{2}x_{21}^{2}x_{22}x_{31}^{2},
x_{11}^{2}x_{12}^{2}x_{21}^{2}x_{22}x_{31}x_{32},
x_{11}^{2}x_{12}^{2}x_{21}^{2}x_{22}x_{32}^{2},
\\
&&
x_{11}^{2}x_{12}^{2}x_{21}x_{22}^{2}x_{31}^{2},
x_{11}^{2}x_{12}^{2}x_{21}x_{22}^{2}x_{31}x_{32},
x_{11}^{2}x_{12}^{2}x_{21}x_{22}^{2}x_{32}^{2},
\\
&&
x_{11}^{2}x_{12}^{2}x_{21}^{2}x_{22}^{2}x_{31},
x_{11}^{2}x_{12}^{2}x_{21}^{2}x_{22}^{2}x_{32}).
\end{eqnarray*}
}
\end{Example}

Now we are interested to study the Cohen-Macaulay property for strong quasi-n-partite graphs.

Let $K$ be a field. We say a graph $G$ is {\em Cohen-Macaulay} over $K$ if $S/I(G)$ has this property.

\begin{Theorem}
\label{Cm}
A quasi-n-partite graph with a loop in each vertex is Cohen-Macaulay.
\end{Theorem}

\begin{proof}
Let $I(\mathcal{K}'_{m_1,\dots,m_n})$ be the edge ideal of a strong quasi-n-partite graph $\mathcal{K}'_{m_1,\dots,m_n}$.
Let $G(I(\mathcal{K}'_{m_1,\dots,m_n}))$ be the unique minimal set of monomial generators of the edge ideal $I(\mathcal{K}'_{m_1,\dots,m_n})$.

A vertex cover of $I(\mathcal{K}'_{m_1,\dots,m_n})$ is a subset $C$ of $\{x_{11},\dots,x_{1m_{1}},\dots,x_{n1},\dots,x_{nm_{n}}\}$ such that
each $u\in G(I(\mathcal{K}'_{m_1,\dots,m_n}))$ is divided by some $x_{ij}\in C$.
A vertex cover $C$ is called {\em minimal} if $C$ is a
vertex of $\mathcal{K}'_{m_1,\dots,m_n}$, and no proper subset of $C$ is a vertex cover of $\mathcal{K}'_{m_1,\dots,m_n}$.
Then the minimal cardinality of the vertex covers of $I(\mathcal{K}'_{m_1,\dots,m_n})$ is
\[
\height(I(\mathcal{K}'_{m_1,\dots,m_n}))=m_{1}+\dots+m_{n}.
\]
Therefore
\[
\depth(T/I(\mathcal{K}'_{m_1,\dots,m_n}))=\dim (T/I(\mathcal{K}'_{m_1,\dots,m_n})),
\]
hence $T/I(\mathcal{K}'_{m_1,\dots,m_n}))$ is a Cohen-Macaulay graph.
\end{proof}

Let $m_{1},\ldots,m_{n}$ be integers, and let
\[
[m_{1}+\dots+m_{n}]=\{1,2,\ldots,m_{1}+\dots+m_{n}\}.
\]
Let $\mathbb{Z}_{+}$ be the set of nonnegative integers and
$\mathbb{Z}^{m_1+\dots+m_{n}}_{+}$ the set of the vectors
$(a_1;a_2;\ldots;a_n)$ with $a_1\in\mathbb{Z}^{m_1}_{+},a_2\in\mathbb{Z}^{m_2}_{+},\ldots,a_n\in\mathbb{Z}^{m_n}_{+}$,
i.e.,
\[
(a_1;a_2;\ldots;a_n)=(a_{11},\ldots,a_{1m_1};a_{21},\ldots,a_{2m_2};\ldots;a_{n1},\ldots,a_{nm_n}),
\]
with $a_{ij}\geq0$ for all $i=1,\ldots,n$ and $j=1,\ldots,m_i$.
Let $(a_1;\ldots;a_n)$ and $(b_1;\ldots;b_n)$ be two vectors belonging to
$\mathbb{Z}^{m_1+\dots+m_{n}}_{+}$. Then we write $(a_1;\ldots;a_n) \geq (b_1;\ldots;b_n)$
if all components $(a_{11}-b_{11},\ldots,a_{nm_n}-b_{nm_{n}})$
of the vector $(a_{1}-b_{1};\ldots;a_{n}-b_{n})$ are nonnegative. Moreover, we write $(a_1;\ldots;a_n) > (b_1;\ldots;b_n)$ if
\[
(a_1;\ldots;a_n) \geq (b_1;\ldots;b_n)
\]
and $(a_1;\ldots;a_n) \neq (b_1;\ldots;b_n)$. We say that $(b_1;\ldots;b_n)$ is a {\em subvector} of $(a_1;\ldots;a_n)$
if $(a_1;\ldots;a_n) \geq (b_1;\ldots;b_n)$.
We set
\[(a_1;\ldots;a_n) \vee (b_1;\ldots;b_n) =(\max \{a_{11},b_{11}\},\ldots,\max\{a_{nm_n},b_{nm_n}\}).
\]
Hence $(a_1;\ldots;a_n) \leq (a_1;\ldots;a_n) \vee (b_1;\ldots;b_n)$ and
\[
(b_1;\ldots;b_n)\leq(a_1;\ldots;a_n) \vee (b_1;\ldots;b_n).
\]
The {\em modulus} of $(a_1;\ldots;a_n) \in \mathbb{Z}^{m_1+\dots+m_{n}}_{+}$
is the number
\[
|(a_1;\ldots;a_n)|=|a_1| +\dots+|a_n|=\sum_{j=1}^{m_{1}}a_{1j}+\dots+\sum_{l=1}^{m_n}a_{nl}.
\]

In \cite{LM}, the following notion of a generalized discrete bi-polymatroid was introduced as a generalization of the discrete bi-polymatroid.

\begin{Definition}
\label{discrete generalized bi-polymatroid}
A {\em generalized discrete bi-polymatroid} on the set $[m_{1}+\dots+m_{n}]$ is a nonempty finite subset
$P^{*} \subset \mathbb{Z}^{m_1+\dots+m_{n}}_{+}$
satisfying the following conditions:

(D1) together with each $(a_{1};\ldots; a_{n}) \in P^{*}$, the set $P^{*}$ contains all its integral subvectors, that is,
if $(a_{1};\ldots; a_{n}) \in P^{*}$ and $(b_{1};\ldots;b_{n}) \in \mathbb{Z}^{m_1+\dots+m_{n}}_{+}$
with $(b_{1};\ldots;b_{n})\leq (a_{1};\ldots;a_{n})$, then $(b_{1};\ldots;b_{n}) \in P^{*}$;

(D2) for all $(a_{1};\ldots; a_{n}) , (b_{1};\ldots;b_{n})\in P^{*}$  with $|(a_{1};\ldots; a_{n})| < |(b_{1};\ldots;b_{n})|$,
there is a vector $(u_{1};\ldots;u_{n}) \in P^{*}$ such that
\[
(a_{1};\ldots; a_{n}) < (u_{1};\ldots; u_{n}) < (a_{1};\ldots; a_{n}) \vee (b_{1};\ldots;b_{n}).
\]
\end{Definition}

A {\em base} of a generalized discrete bi-polymatroid $P^{*}$ is a vector $(a_{1};\ldots; a_{n})\in P^{*}$
such that $(a_{1};\ldots; a_{n}) < (b_{1};\ldots;b_{n})$
for no $(b_{1};\ldots;b_{n})\in P^{*}$. The set of bases of $P^{*}$ is denoted by $B(P^{*})$.
It follows from (D2) that if $(a_{1};\ldots; a_{n})$ and $(b_{1};\ldots;b_{n})$ are bases of $P^{*}$, then
$|(a_{1};\ldots; a_{n})| = |(b_{1};\ldots; b_{n})|$. The nonnegative integer
$|(a_{1};\ldots; a_{n})|$ with $(a_{1};\ldots; a_{n}) \in B(P^{*})$ is called the {\em rank} of $P^{*}$.

\medskip

For the canonical basis vector $\varepsilon _{i}\in \mathbb{R}^{n}$, we set $\varepsilon_{i}^{*}=\{\varepsilon_{i1},\ldots,\varepsilon_{im_i}\}$
which $\varepsilon_{ij}$ is a canonical basis vector of $\mathbb{R}^{m_1+\dots+m_{n}}_{+}$
with the $(m_1+\dots+m_{i-1}+j)$-th component is equal to 1 and zero for other components.

A characterization of generalized discrete bi-polymatroids in terms of their bases is the following:

\begin{Theorem}
\label{generalized bi-excheng property}
\cite[Generalized bi-exchange property]{LM}

$P^{*} \subset \mathbb{Z}^{m_1+\dots+m_{n}}_{+}$ is a generalized discrete  bi-polymatroid if and only if
for any $(a_{1};\ldots ;a_{n})$, $(b_{1};\ldots;b_{n})\in  B(P^{*})$ with $a_{ij} > b_{ij}$, then there exists $l\in \{1,\dots,m_{i}\}$
with $a_{il} < b_{il}$ such that $(a_{1};\ldots; a_{n})-\varepsilon_{ij}+\varepsilon_{il} \in B(P^{*})$.
\end{Theorem}

A monomial ideal generated by all the monomials corresponding to the set $B(P^{*})$
of bases of a generalized discrete bi-polymatroid is called {\em generalized bi-polymatroidal ideal}
and it is generated by all the monomials $\underline{\Xb}^{(a_1;\dots;a_n)}$ with $(a_1;\dots;a_n) \in B(P^{*})$,
where $\underline{\Xb}^{(a_1;\dots;a_n)}$ stands for
$x_{11}^{a_{11}} \cdots x_{1m_1}^{a_{1m_{1}}} \cdots x_{n1}^{a_{n1}} \cdots x_{nm_n}^{a_{nm_{n}}}$.

In \cite{LM}, the following definition of generalized bi-polymatroidal ideal as a consequence of Theorem \ref{generalized bi-excheng property}.

\begin{Definition}
\label{generalized bi-polymatroidal}
A monomial ideal $L^{*}$ of $T$ generated in a single degree is said to be {\em generalized bi-polymatroidal}
if the following condition is satisfied: for all monomials
$u = x_{11}^{a_{11}} \cdots x_{1m_1}^{a_{1m_{1}}} \cdots x_{n1}^{a_{n1}} \cdots x_{nm_n}^{a_{nm_{n}}}$
and
$v = x_{11}^{b_{11}} \cdots x_{1m_1}^{b_{1m_{1}}} \cdots x_{n1}^{b_{n1}} \cdots x_{nm_n}^{b_{nm_{n}}}$
in $G(L^{*})$ and for each $ij$ with $a_{ij} > b_{ij}$, then there exists $l\in \{1,\dots,m_{i}\}$ with $a_{il}<b_{il}$
such that $x_{il}(u/x_{ij})\in G(L^{*})$.
\end{Definition}

Next we want to study powers of generalized graph ideals of strong quasi-n-partite graphs.

\begin{Proposition}
\label{power}
Let $T=K[x_{11},\ldots,x_{1m_1},x_{21},\ldots,x_{2m_2},\ldots,x_{n1},\ldots,x_{nm_n}]$, and let
\[
I_{t}(\mathcal{K}'_{m_1,\dots,m_n})=\sum_{h=1}^{r}\prod _{i=1}^{n} L_{i,q_{h_{i}},2}, \quad q_{h_{i}}\geq 0
\]
be the generalized graph ideal of a strong quasi-n-partite graph $\mathcal{K}'_{m_1,\dots,m_n}$,
where the ideals $L_{i,q_{h_{i}},2}$ are Veronese type ideals of degree $q_{h_{i}}$ generated by the monomials $x_{i1}^{a_{i1}} \ldots x_{im_{i}}^{a_{im_{i}}}$ with $\sum_{j=1}^{m_{i}}a_{ij}=q_{h_{i}}$ and $0 \leq a_{ij}\leq 2$ for $i=1,\dots,n$.
Suppose that $\sum_{i=1}^{n}q_{h_{i}}=t$ and $q_{h_{i}}\neq t$ for all $i=1,\dots,n$ and $h=1,\dots,r$.
Then for all $k\geq 1$ we have
\[
I_{t}(\mathcal{K}'_{m_1,\dots,m_n})^{k}=\sum_{k_1,\dots,k_r\geq 0,\sum_{h=1}^r k_h=k}\prod _{i=1}^n\prod_{h=1}^rL_{i,q_{h_{i}},2}^{k_{h}},
\]
where $0\leq q_{h_{i}}\leq 2m_{i}$ for $i=1,\dots,n$.
\end{Proposition}

\begin{proof}
Let $\mathcal{K}'_{m_1,\dots,m_n}$ be a strong quasi-n-partite graph with vertex set $V=V_1\union V_2\union \cdots\union V_n$ and $V_i=\{x_{i1},\ldots,x_{im_i}\}$ for $i=1,\ldots,n$.

Suppose that $I_{t}(\mathcal{K}'_{m_1,\dots,m_n})=\sum_{h=1}^{r}L_{1,q_{h_{1}},2}\cdots L_{n,q_{h_{n}} ,2}$
be the generalized graph ideal of $\mathcal{K}'_{m_1,\dots,m_n}$, where $\sum_{i=1}^{n}q_{h_{i}}=t$ and $q_{h_{i}}\neq t$ for $h=1,\dots,r$ .
Then for all $k\geq 1$ we have
\begin{eqnarray*}
I_{t}(\mathcal{K}'_{m_1,\dots,m_n})^{k}
&=& \sum_{k_1,\dots,k_r\geq 0,\sum _{h=1}^r k_h=k}(L_{1,q_{1_{1}},2}^{k_1}L_{2,q_{1_{2}},2}^{k_1}
\cdots L_{n,q_{1_{n}},2}^{k_1})\cdots
\\
&&
(L_{1,q_{{r-1}_{1}},2}^{k_{r-1}}L_{2,q_{{r-1}_{2}},2}^{k_{r-1}}\cdots L_{n,q_{{r-1}_{n}},2}^{k_{r-1}})(L_{1,q_{r_{1}},2}^{k_r}L_{2,q_{r_{2}},2}^{k_r}\cdots L_{n,q_{r_{n}},2}^{k_r})\\
&=& \sum_{k_1,\dots,k_r\geq 0,\sum _{h=1}^r k_h=k}(L_{1,q_{1_{1}},2}^{k_1}L_{1,q_{2_{1}},2}^{k_2}
\cdots L_{1,q_{r_{1}},2}^{k_r})\cdots
\\
&&
(L_{n-1,q_{1_{n-1}},2}^{k_1}L_{n-1,q_{2_{n-1}},2}^{k_2}\cdots L_{n-1,q_{r_{n-1}},2}^{k_{r}})(L_{n,q_{1_{n}},2}^{k_1}L_{n,q_{2_{n}},2}^{k_2}\cdots L_{n,q_{r_{n}},2}^{k_{r}})\\
&=& \sum_{k_1,\dots,k_r\geq 0,\sum_{h=1}^r k_h=k}\prod _{i=1}^n\prod_{h=1}^rL_{i,q_{h_{i}},2}^{k_{h}}.
\end{eqnarray*}
Thus the desired conclusion follows.
\end{proof}

\begin{Lemma}
\label{product}
Let $I_{t}(\mathcal{K}'_{m_1,\dots,m_n})$ and $I_{t'}(\mathcal{K}'_{m_1,\dots,m_n})$ be two generalized graph ideals associated to $\mathcal{K}'_{m_1,\dots,m_n}$. Then $I_{t}(\mathcal{K}'_{m_1,\dots,m_n})I_{t'}(\mathcal{K}'_{m_1,\dots,m_n})$ is a generalized bi-polymatroidal ideal.
\end{Lemma}

\begin{proof}
Let $\mathcal{K}'_{m_1,\dots,m_n}$ be a strong quasi-n-partite graph with vertex set $V=V_1\union V_2\union \cdots\union V_n$ and $V_i=\{x_{i1},\ldots,x_{im_i}\}$ for $i=1,\ldots,n$.
Let $t, q_{1}, \dots ,q_{n}$ be non negative integers and $\sum_{i=1}^{n}q_{i}=t$, $q_{1}, \dots ,q_{n}\geq 0$. Assume further that $q_{i}\neq t$ for $i=1,\dots,n$.
It then follows that
\[
B_{t,2}=\{(a_{1};\dots;a_{n})\in \ZZ^{m_{1}+\dots+m_{n}}_{+}:|a_{i}|=\sum _{j=1}^{m_{i}}a_{ij}=q_{i}, \quad 0\leq a_{ij}\leq 2\}
\]
is the set of bases of a generalized discrete bi-polymatroid of rank $t$.
In fact, let $(a_{1};\ldots ;a_{n})$, $(b_{1};\ldots;b_{n})\in  B_{t,2}$ with $a_{ij} > b_{ij}$, for some $l\in \{1,\dots,m_{i}\}$
with $a_{il} < b_{il}$, one has $(a_{1};\ldots; a_{n})-\varepsilon_{ij}+\varepsilon_{il} \in B_{t,2}$.

The generalized bi-polymatroidal ideal corresponding to the set $B_{t,2}$ is
\[
I_{t}(\mathcal{K}'_{m_1,\dots,m_n})=\sum_{\sum_{i=1}^{n}q_{i}=t, q_{i}\neq t}L_{1,q_{1},2}\cdots L_{n,q_{n},2}\subset T,
\]
where the ideals $L_{i,q_{i},2}$ are Veronese type ideals of degree $q_{i}$ generated by the monomials $x_{i1}^{a_{i1}} \cdots x_{im_{i}}^{a_{im_{i}}}$
with $\sum_{j=1}^{m_{i}}a_{ij}=q_{i}$ and $0 \leq a_{ij}\leq 2$ for all $i=1,\dots,n$. Here we set $L_{i,q_{i},2}=(0)$, if $q_{i}=t$.

Let $\alpha$ and $\beta$ be two monomials in the same degree.
We set
\[
d(\alpha,\beta)=\frac{1}{2}\sum_{i}|\deg_{x_{i}}(\alpha)-\deg_{x_{i}}(\beta)|.
\]
Hence, $d(\alpha,\beta)=0$ if and only if $\alpha=\beta$.

Let $\alpha,\alpha_{1}\in G(I_{t}(\mathcal{K}'_{m_1,\dots,m_n}))$ and $\beta,\beta_{1}\in G(I_{t'}(\mathcal{K}'_{m_1,\dots,m_n}))$ and suppose that
\[
\deg_{x_{ij}}(\alpha_{1}\beta_{1})> \deg_{x_{ij}}(\alpha \beta).
\]
It may be assumed that $\deg_{x_{ij}}(\alpha_{1})>\deg_{x_{ij}}(\alpha)$. Since $I_{t}(\mathcal{K}'_{m_1,\dots,m_n})$ is a generalized bi-polymatroidal ideal, it follows that there exists an integer $il_{1}$
such that $\deg_{x_{il_{1}}}(\alpha)>\deg_{x_{il_{1}}}(\alpha_{1})$ and
\[
\alpha_{2}=x_{il_{1}}(\alpha_{1}/x_{ij})\in G(I_{t}(\mathcal{K}'_{m_1,\dots,m_n})).
\]
Then we have $d(\alpha_{2},\alpha)<d(\alpha,\alpha_{1})$.

If $\deg_{x_{il_{1}}}(\beta) \geq \deg_{x_{il_{1}}}(\beta_{1})$ we are done, because then $\deg_{x_{il_{1}}}(\alpha \beta)> \deg_{x_{il_{1}}}(\alpha_{1}\beta_{1})$, and
\[
x_{il_{1}}(\alpha_{1}\beta_{1}/x_{ij})=\alpha_{2}\beta_{1}\in G(I_{t}(\mathcal{K}'_{m_1,\dots,m_n})I_{t'}(\mathcal{K}'_{m_1,\dots,m_n})).
\]
Otherwise, $\deg_{x_{il_{1}}}(\beta_{1})> \deg_{x_{il_{1}}}(\beta)$. Thus there exists $it_{1}$ with $\deg_{x_{it_{1}}}(\beta)> \deg_{x_{it_{1}}}(\beta_{1})$ and such that
\[
\beta_{2}=x_{it_{1}}(\beta_{1}/x_{il_{1}})\in G(I_{t'}(\mathcal{K}'_{m_1,\dots,m_n})).
\]
Thus $d(\beta_{2},\beta)<d(\beta_{1},\beta)$.

If $\deg_{x_{it_{1}}}(\alpha)\geq \deg_{x_{it_{1}}}(\alpha_{2})$, then $\deg_{x_{it_{1}}}(\alpha \beta)>\deg_{x_{it_{1}}}(\alpha_{2}\beta_{1})=\deg_{x_{it_{1}}}(x_{il_{1}}(\alpha_{1}\beta_{1}/x_{ij}))$. Hence,
\[
\deg_{x_{it_{1}}}(\alpha \beta) > \deg_{x_{it_{1}}}(\alpha_{1}\beta_{1})
\]
if $ij\neq it_{1}$. Then
\[
x_{it_{1}}(\alpha_{1}\beta_{1}/x_{ij})=\alpha_{2}\beta_{2}\in G(I_{t}(\mathcal{K}'_{m_1,\dots,m_n})I_{t'}(\mathcal{K}'_{m_1,\dots,m_n})).
\]
On the other hand, if $it_{1}=ij$, then $\alpha_{1}\beta_{1}=\alpha_{2}\beta_{2}$, and by induction we suppose that the generalized bi-exchange property holds since
$d(\alpha_{2},\alpha)<d(\alpha_{1},\alpha)$ and $d(\beta_{2},\beta)< d(\beta_{1},\beta)$.

Otherwise $\deg_{x_{it_{1}}}(\alpha_{2})> \deg_{x_{it_{1}}}(\alpha)$. Then there exists $il_{2}$ with $\deg_{x_{il_{2}}}(\alpha) > \deg_{x_{il_{2}}}(\alpha_{2})$ and such that
\[
\alpha_{3}=x_{il_{2}}(\alpha_{2}/x_{it_{1}})\in G(I_{t}(\mathcal{K}'_{m_1,\dots,m_n})).
 \]
If $\deg_{x_{il_{2}}}(\beta)\geq \deg_{x_{il_{2}}}(\beta_{2})$, hence $\deg_{x_{il_{2}}}(\alpha \beta)> \deg_{x_{il_{2}}}(\alpha_{2}\beta_{2})=\deg_{x_{il_{2}}}(x_{it_{1}}(\alpha_{1}\beta_{1}/x_{ij}))$. Then if $il_{2}\neq ij$, thus
\[
\deg_{x_{il_{2}}}(\alpha \beta)> \deg_{x_{il_{2}}}(\alpha_{1}\beta_{1}),
\]
and we are done since
$x_{il_{2}}(\alpha_{1}\beta_{1}/x_{ij})=\alpha_{3}\beta_{2}\in I_{t}(\mathcal{K}'_{m_1,\dots,m_n})I_{t'}(\mathcal{K}'_{m_1,\dots,m_n})$.

On the other hand if $il_{2}=ij$, hence $\alpha_{3}\beta_{2}=\alpha_{1}\beta_{1}$, and by induction on the distance we have the desired generalized bi-exchange property.
Otherwise $\deg_{x_{il_{2}}}(\beta_{2})> \deg_{x_{il_{2}}}(\beta)$.

We proceed in this way. Suppose we have constructed sequences $x_{il_{1}},\dots, x_{il_{z}}$, $x_{it_{1}},\dots, x_{it_{z-1}}$, and
$\alpha_{1},\dots,\alpha_{z+1}\in G(I_{t}(\mathcal{K}'_{m_1,\dots,m_n}))$ and $\beta_{1},\dots,\beta_{z}\in G(I_{t'}(\mathcal{K}'_{m_1,\dots,m_n}))$ such that for $q=1,\dots,z$ we have
\begin{enumerate}
\item[(a)] $x_{it_{q-1}}$ divides $\alpha_{q}$ and $x_{il_{q}}$ divides $\beta_{q}$,

\item[(b)] $\alpha_{q+1}=x_{il_{q}}(\alpha_{q}/x_{it_{q-1}})$ and $\beta_{q}=x_{it_{q-1}}(\beta_{q-1}/x_{il_{q-1}})$,

\item[(c)] $d(\alpha_{q+1}, \alpha)<d(\alpha_{q},\alpha)$ for $q\neq z$, $d(\beta_{q+1},\beta)< d(\beta_{q},\beta)$,

\item[(d)] $\deg_{x_{il_{q}}}(\alpha)> \deg_{x_{il_{q}}}(\alpha_{q})$ and $\deg_{x_{it_{q}}}(\beta)> \deg_{x_{it_{q}}}(\beta_{q})$.

\end{enumerate}
Furthermore we set $it_{0}=ij$ for symmetric reason.
Then
\[
\alpha_{z+1}=x_{il_{z}}\cdots x_{il_{1}}(\alpha_{1}/x_{ij}x_{it_{1}}\cdots x_{it_{z-1}})\quad  \text{and} \quad  \beta_{z}=x_{it_{z-1}}\cdots x_{it_{1}}(\beta_{1}/x_{il_{1}}\cdots x_{il_{z-1}}).
\]
If $\deg_{x_{il_{z}}}(\beta)\geq \deg_{x_{il_{z}}}(\beta_{z})$, thus $\deg_{x_{il_{z}}}(\alpha \beta) > \deg_{x_{il_{z}}}(\alpha_{z}\beta_{z})=\deg_{x_{il_{z}}}(x_{it_{z-1}}(\alpha_{1}\beta_{1}/x_{ij}))$ by (d).
Hence, if $il_{z}\neq ij$, then $\deg_{x_{il_{z}}}(\alpha \beta) > \deg_{x_{il_{z}}}(\alpha_{1}\beta_{1})$, and we are done since
\[
x_{il_{z}}(\alpha_{1}\beta_{1}/x_{ij})=\alpha_{z+1}\beta_{z}\in G(I_{t}(\mathcal{K}'_{m_1,\dots,m_n})I_{t'}(\mathcal{K}'_{m_1,\dots,m_n})).
\]
On the other hand, if $il_{z}=ij$, and hence $\alpha_{1}\beta_{1}=\alpha_{z+1}\beta_{z}$ and by induction on the distance we have the desired generalized bi-exchange property.

Otherwise $\deg_{x_{il_{z}}}(\beta_{z})> \deg_{x_{il_{z}}}(\beta)$, and there exists $it_{z}$ with $\deg_{x_{it_{z}}}(\beta)> \deg_{x_{it_{z}}}(\beta_{z})$ and such that
$\beta_{z+1}=x_{it_{z}}(\beta_{z}x_{il_{z}})\in G(I_{t'}(\mathcal{K}'_{m_1,\dots,m_n}))$. In addition, we have $d(\beta_{z+1},\beta)< d(\beta_{z},\beta)$. Therefore, the new elements $x_{it_{z}}$ and $\beta_{z+1}$ satisfy
the properties (a)-(d).

If $\deg_{x_{it_{z}}}(\alpha)> \deg_{x_{it_{z}}}(\alpha_{z+1})$, then by (d),
\[
\deg_{x_{it_{z}}}(\alpha \beta) > \deg_{x_{it_{z}}}(\alpha_{z+1}\beta_{z})= \deg_{x_{it_{z}}}(x_{il_{z}}(\alpha_{1}\beta_{1}/x_{ij})).
\]
Hence, if $it_{z}\neq ij$, thus $\deg_{x_{it_{z}}}(\alpha \beta)> \deg_{x_{it_{z}}}(\alpha_{1}\beta_{1})$, and we are done since
\[
x_{it_{z}}(\alpha_{1}\beta_{1}/x_{ij})=\alpha_{z+1}\beta_{z+1}\in G(I_{t}(\mathcal{K}'_{m_1,\dots,m_n})I_{t'}(\mathcal{K}'_{m_1,\dots,m_n})).
\]
On the other hand, if $it_{z}=ij$, and then $\alpha_{1}\beta_{1}=\alpha_{z+1}\beta_{z+1}$ and by induction on the distance we have the desired generalized bi-exchange property.

Otherwise, $\deg_{x_{it_{z}}}(\alpha_{z+1}) > \deg_{x_{it_{z}}}(\alpha)$, and there exists $il_{z+1}$ with $\deg_{x_{il_{z}}}(\alpha)> \deg_{x_{il_{z}}}(\alpha_{z+1})$ and such that
$\alpha_{z+2}=x_{il_{z+1}}(\alpha_{z+1}/x_{it_{z}})\in G(I_{t}(\mathcal{K}'_{m_1,\dots,m_n}))$.

Moreover, $d(\alpha_{z+2},\alpha)< d(\alpha_{z+1},\alpha)$. Then we have the conditions (a)-(d) as before but $z$
replace by $z+1$. By condition (c) we conclude that the process must be terminated.
\end{proof}

As a consequence of Lemma \ref{product} we have the following corollary.

\begin{Corollary}
\label{power-bi-polymatroid}
Let $T=K[x_{11},\ldots,x_{1m_{1}},\ldots,x_{n1},\ldots,x_{nm_n}]$, and let $I_{t}(\mathcal{K}'_{m_1,\dots,m_n})\subset T$ be the generalized graph ideal of a strong quasi-n-partite graph $\mathcal{K}'_{m_1,\dots,m_n}$. Then
$I_{t}(\mathcal{K}'_{m_1,\dots,m_n})^{k}$ is a generalized bi-polymatroidal ideal for all $k\geq 1$.
\end{Corollary}

\begin{Example}
\label{quasi}
{
Let $T=K[x_{11},x_{12},x_{21},x_{22},x_{31},x_{32}]$ be a polynomial ring and $\mathcal{K}'_{2,2,2}$ be the strong quasi-3-partite graph on vertices $x_{11},x_{12},x_{21},x_{22},x_{31},x_{32}$. Then
\begin{eqnarray*}
I_{11}(\mathcal{K}'_{2,2,2})
&=&L_{1,3,2}L_{2,4,2}L_{3,4,2}+L_{1,4,2}L_{2,3,2}L_{3,4,2}+L_{1,4,2}L_{2,4,2}L_{3,3,2}\\
&=&(x_{11}^{2}x_{12}x_{21}^{2}x_{22}^{2}x_{31}^{2}x_{32}^{2},
x_{11}x_{12}^{2}x_{21}^{2}x_{22}^{2}x_{31}^{2}x_{32}^{2},
x_{11}^{2}x_{12}^{2}x_{21}^{2}x_{22}x_{31}^{2}x_{32}^{2},\\
&&
x_{11}^{2}x_{12}^{2}x_{21}x_{22}^{2}x_{31}^{2}x_{32}^{2},
x_{11}^{2}x_{12}^{2}x_{21}^{2}x_{22}^{2}x_{31}^{2}x_{32},
x_{11}^{2}x_{12}^{2}x_{21}^{2}x_{22}^{2}x_{31}x_{32}^{2})
\end{eqnarray*}
is a generalized bi-polymatroidal ideal.
Thus by using Proposition \ref{power} it follows that
\begin{eqnarray*}
I_{11}(\mathcal{K}'_{2,2,2})^{2}
&=&L_{1,6,4}L_{2,8,4}L_{3,8,4}+L_{1,8,4}L_{2,6,4}L_{3,8,4}+L_{1,8,4}L_{2,8,4}L_{3,6,4}\\
&&
+L_{1,7,4}L_{2,7,4}L_{3,8,4}+L_{1,8,4}L_{2,7,4}L_{3,7,4}+L_{1,7,4}L_{2,8,4}L_{3,7,4}.
\end{eqnarray*}
Therefore, Corollary \ref{power-bi-polymatroid} implies that
\begin{eqnarray*}
I_{11}(\mathcal{K}'_{2,2,2})^{2}
&=&(x_{11}^{4}x_{12}^{2}x_{21}^{4}x_{22}^{4}x_{31}^{4}x_{32}^{4},
x_{11}^{2}x_{12}^{4}x_{21}^{4}x_{22}^{4}x_{31}^{4}x_{32}^{4},
x_{11}^{4}x_{12}^{4}x_{21}^{4}x_{22}^{2}x_{31}^{4}x_{32}^{4},\\
&&
x_{11}^{4}x_{12}^{4}x_{21}^{2}x_{22}^{4}x_{31}^{4}x_{32}^{4},
x_{11}^{4}x_{12}^{4}x_{21}^{4}x_{22}^{4}x_{31}^{4}x_{32}^{2},
x_{11}^{4}x_{12}^{4}x_{21}^{4}x_{22}^{4}x_{31}^{2}x_{32}^{4}\\
&&
x_{11}^{3}x_{12}^{3}x_{21}^{4}x_{22}^{4}x_{31}^{4}x_{32}^{4},
x_{11}^{4}x_{12}^{3}x_{21}^{4}x_{22}^{3}x_{31}^{4}x_{32}^{4},
x_{11}^{4}x_{12}^{3}x_{21}^{3}x_{22}^{4}x_{31}^{4}x_{32}^{4},\\
&&
x_{11}^{4}x_{12}^{3}x_{21}^{4}x_{22}^{4}x_{31}^{4}x_{32}^{3},
x_{11}^{4}x_{12}^{3}x_{21}^{4}x_{22}^{4}x_{31}^{3}x_{32}^{4},
x_{11}^{3}x_{12}^{4}x_{21}^{4}x_{22}^{3}x_{31}^{4}x_{32}^{4},\\
&&
x_{11}^{3}x_{12}^{4}x_{21}^{3}x_{22}^{4}x_{31}^{4}x_{32}^{4},
x_{11}^{3}x_{12}^{4}x_{21}^{4}x_{22}^{4}x_{31}^{4}x_{32}^{3},
x_{11}^{3}x_{12}^{4}x_{21}^{4}x_{22}^{4}x_{31}^{3}x_{32}^{4},\\
&&
x_{11}^{4}x_{12}^{4}x_{21}^{3}x_{22}^{3}x_{31}^{4}x_{32}^{4},
x_{11}^{4}x_{12}^{4}x_{21}^{4}x_{22}^{3}x_{31}^{4}x_{32}^{3},
x_{11}^{4}x_{12}^{4}x_{21}^{4}x_{22}^{3}x_{31}^{3}x_{32}^{4},\\
&&
x_{11}^{4}x_{12}^{4}x_{21}^{3}x_{22}^{4}x_{31}^{4}x_{32}^{3},
x_{11}^{4}x_{12}^{4}x_{21}^{3}x_{22}^{4}x_{31}^{3}x_{32}^{4},
x_{11}^{4}x_{12}^{4}x_{21}^{4}x_{22}^{4}x_{31}^{3}x_{32}^{3})
\end{eqnarray*}
is again a generalized bi-polymatroidal ideal.
}
\end{Example}

\section{Monomial localizations of generalized graph ideals}
\label{two}
The main goal of this section is to study some permanence properties of the generalized graph ideals.

Let $K$ be a field, $S=K[x_{1},\dots,x_{n}]$ the polynomial ring in the indeterminates $x_{1},\dots,x_{n}$ and $I\subset S$ a monomial ideal.

Furthermore let $T=K[x_{11},\ldots,x_{1m_1},x_{21},\ldots,x_{2m_2},\ldots,x_{n1},\ldots,x_{nm_n}]$ denote the polynomial ring in $m_{1}+\dots+m_{n}$ variables over $K$.
We first show

\begin{Theorem}
\label{colon-generalized}
Let $I_{t}(\mathcal{K}'_{m_1,\dots,m_n})$ be the generalized graph ideal of a strong quasi-n-partite graph $\mathcal{K}'_{m_1,\dots,m_n}$.
Then $I_{t}(\mathcal{K}'_{m_1,\dots,m_n}):u$ is a generalized bi-polymatroidal ideal for all monomials $u$.
\end{Theorem}

\begin{proof}
Let $\mathcal{K}'_{m_1,\dots,m_n}$ be a strong quasi-n-partite graph with vertex set $V=V_1\union V_2\union \cdots\union V_n$ and $V_i=\{x_{i1},\ldots,x_{im_i}\}$ for $i=1,\ldots,n$.
Assume further that $q_{1}, \dots ,q_{n},t$ be non negative integers and $\sum_{i=1}^{n}q_{i}=t$, $q_{1}, \dots ,q_{n}\geq 0$.
Let
\[
I_{t}(\mathcal{K}'_{m_1,\dots,m_n})=\sum_{\sum_{i=1}^{n}q_{i}=t,q_{i}\neq t}L_{1,q_{1},2}\cdots L_{n,q_{n},2}
 \]
be the generalized graph ideal of a strong quasi-n-partite graph $\mathcal{K}'_{m_1,\dots,m_n}$,
where the ideals $L_{i,q_{i},2}$ are Veronese type ideals of degree $q_{i}$ generated by the monomials $x_{i1}^{a_{i1}} \ldots x_{im_{i}}^{a_{im_{i}}}$ with $\sum_{j=1}^{m_{i}}a_{ij}=q_{i}$ and $0 \leq a_{ij}\leq 2$ for $i=1,\dots,n$.

In fact, let
$u = x_{11}^{a_{11}} \cdots x_{1m_1}^{a_{1m_{1}}} \cdots x_{n1}^{a_{n1}} \cdots x_{nm_n}^{a_{nm_{n}}}$
and
$v = x_{11}^{b_{11}} \cdots x_{1m_1}^{b_{1m_{1}}} \cdots x_{n1}^{b_{n1}} \cdots x_{nm_n}^{b_{nm_{n}}}$
in $G(I_{t}(\mathcal{K}'_{m_1,\dots,m_n}))$ with $a_{ij} > b_{ij}$. Then for some $l\in \{1,\dots,m_{i}\}$ such that $a_{il}<b_{il}$, we have
$x_{il}(u/x_{ij})\in G(I_{t}(\mathcal{K}'_{m_1,\dots,m_n}))$.
We will show that for all
variables $x_{ij}$, $I_{t}(\mathcal{K}'_{m_1,\dots,m_n}) : x_{ij}$ is a generalized bi-polymatroidal ideal
of $T$. In addition, we set
\[
\mathcal{I}_{t}(\mathcal{K}'_{m_1,\dots,m_n}) = I_{t}(\mathcal{K}'_{m_1,\dots,m_n}) : x_{ij}.
\]
Let $u, v \in  G(\mathcal{I}_{t}(\mathcal{K}'_{m_1,\dots,m_n}))$. Thus
\[
x_{ij}u,x_{ij}v \in  x_{ij}\mathcal{I}_{t}(\mathcal{K}'_{m_1,\dots,m_n})\subseteq I_{t}(\mathcal{K}'_{m_1,\dots,m_n}).
\]
If $\deg_{x_{ij}}(u)=\deg_{x_{ij}}(v)$, then $x_{ij}u$ and $x_{ij}v$
have the property described in Definition \ref{generalized bi-polymatroidal}.
Hence if $\deg_{x_{ij}}(u) > \deg_{x_{ij}}(v)$,
then $x_{ij}$ divides $u$.

For a variable $x_{ir}$ with $\deg_{x_{ir}}(u)>\deg_{x_{ir}}(v)$,
we show that there exists a variable $x_{ir'}$ with $\deg_{x_{ir'}}(v)> \deg_{x_{ir'}}(u)$
and such that
\[
x_{ir'}(u/x_{ir}) \in G(\mathcal{I}_{t}(\mathcal{K}'_{m_1,\dots,m_n})).
\]
We know that $\deg _{x_{ir}} (x_{ij}u) > \deg _{x_{ir}} (x_{ij}v)$ and $I_{t}(\mathcal{K}'_{m_1,\dots,m_n})$ is
a generalized bi-polymatroidal ideal, it then follows that there exists $x_{ir'}$ with $\deg _{x_{ir'}} (x_{ij}u) < \deg _{x_{ir'}} (x_{ij}v)$
such that $x_{ir'}(x_{ij}u/x_{ir}) \in G(I_{t}(\mathcal{K}'_{m_1,\dots,m_n}))$. Then $x_{ir'}(u/x_{ir})\in \mathcal{I}_{t}(\mathcal{K}'_{m_1,\dots,m_n})$. The same result holds
true for $\deg_{x_{ij}} (u) < \deg_{x_{ij}} (v)$.
\end{proof}

Let $I\subset S$ be a monomial ideal.
We say that $I$ has {\em linear quotients}, if there is an ordering $u_1, \ldots, u_r$
of the monomials belonging to $G(I)$ with $\deg(u_1)\leq \cdots \leq \deg(u_r)$ such that for
each $2 \leq j \leq r$, the colon ideal $(u_1,\dots, u_{j-1}) : u_j$ is generated by a subset of
$\{x_{1},\dots,x_{n}\}$.
The reader can find more information in \cite[Definition 6.3.45]{VI}.

\begin{Theorem}
\label{linear}
Let $I_{t}(\mathcal{K}'_{m_1,\dots,m_n})\subset T$ be the generalized graph ideal of a strong quasi-n-partite graph $\mathcal{K}'_{m_1,\dots,m_n}$.
If $I_{t}(\mathcal{K}'_{m_1,\dots,m_n}): u$ is generated in a single degree and has linear quotients for all monomials $u$, then $I_{t}(\mathcal{K}'_{m_1,\dots,m_n})$ is a generalized bi-polymatroidal ideal.
\end{Theorem}

\begin{proof}
Suppose that $I_{t}(\mathcal{K}'_{m_1,\dots,m_n}): u$ is generated in a single degree and has linear quotients for all monomials $u$. By \cite[Lemma 4.1]{CH} we have $I_{t}(\mathcal{K}'_{m_1,\dots,m_n}) : u$ has a linear resolution.
Thus $I_{t}(\mathcal{K}'_{m_1,\dots,m_n}) : u$ is generated in a single degree for all monomials $u$ of $T$.

Let $u',u'' \in G(I_{t}(\mathcal{K}'_{m_1,\dots,m_n}))$ with $\deg_{x_{ij}}(u') > \deg_{x_{ij}} (u'')$.
We will show that there exists a variable $x_{ir}$ such that $\deg_{x_{ir}} (u'') > \deg_{x_{ir}} (u')$ and $x_{ir}(u'/x_{ij})\in G(I_{t}(\mathcal{K}'_{m_1,\dots,m_n}))$.

Now assume in addition that $\deg_{x_{ij}} (u') > \deg_{ x_{ij}} (u'')$. Then $I_{t}(\mathcal{K}'_{m_1,\dots,m_n}) : \frac{u'}{x_{ij}}$
is generated in the same degree.
Since $x_{ij} \in G(I_{t}(\mathcal{K}'_{m_1,\dots,m_n}): \frac{u'}{x_{ij}})$,
it then follows that $I_{t}(\mathcal{K}'_{m_1,\dots,m_n}) : \frac{u'}{x_{ij}}$ is generated in degree 1.
Hence, we know that
\[
u''/\GCD (u'',\frac{u'}{x_{ij}}) \in I_{t}(\mathcal{K}'_{m_1,\dots,m_n}): \frac{u'}{x_{ij}},
\]
there exists $f \in G(I_{t}(\mathcal{K}'_{m_1,\dots,m_n}))$ such that $x_{ir} = f/\GCD (f,\frac{u'}{x_{ij}})$
for some $r$ and $x_{ir}$ divides $u''/\GCD (u'',\frac{u'}{x_{ij}})$.
One has $\deg _{x_{ir}} (u'') > \deg_{x_{ir}}(\frac{u'}{x_{ij}})$.
Since $\deg_{x_{ij}} (u') > \deg_{x_{ij}}  (u'')$, it follows that $x_{ij} \neq x_{ir}$. Hence, $\deg _{x_{ir}} (u'') > \deg_{x_{ir}}(\frac{u'}{x_{ij}})=\deg_{x_{ir}}(u')$.
Our assumption implies that $I_{t}(\mathcal{K}'_{m_1,\dots,m_n})$ is generated in a single degree.
Thus $\deg(f) = \deg(u')$. In fact, since $x_{ir}=f/\GCD (f, \frac{u'}{x_{ij}})$,
this implies that
\[
\deg_{ x_{ir'}}(f) \leq  \deg_{x_{ir'}}(\frac{u'}{x_{ij}})=\deg _{x_{ir'}}(u')
\]
for all $r' \neq r,j$,
$\deg _{x_{ir}}(f) = \deg _{x_{ir}}(\frac{u'}{x_{ij}})+1=\deg _{x_{ir}}(u')+1$,
and
\[
\deg_{x_{ij}}(f) \leq \deg _{x_{ij}}(\frac{u'}{x_{ij}})=\deg_{x_{ij}}(u')-1.
\]
Consequently, we obtain $f=x_{ir}(u'/x_{ij})$, where $x_{ir}(u'/x_{ij}) \in G(I_{t}(\mathcal{K}'_{m_1,\dots,m_n}))$.
\end{proof}

\begin{Example}
\label{Example-colon}
{\em
Let $T=K[x_{11},x_{12},x_{21},x_{22},x_{31},x_{32}]$. Consider the generalized graph ideal of Example \ref{quasi}:
\begin{eqnarray*}
I_{11}(\mathcal{K}'_{2,2,2})
&=&L_{1,3,2}L_{2,4,2}L_{3,4,2}+L_{1,4,2}L_{2,3,2}L_{3,4,2}+L_{1,4,2}L_{2,4,2}L_{3,3,2}\\
&=&(x_{11}^{2}x_{12}x_{21}^{2}x_{22}^{2}x_{31}^{2}x_{32}^{2},
x_{11}x_{12}^{2}x_{21}^{2}x_{22}^{2}x_{31}^{2}x_{32}^{2},
x_{11}^{2}x_{12}^{2}x_{21}^{2}x_{22}x_{31}^{2}x_{32}^{2},\\
&&
x_{11}^{2}x_{12}^{2}x_{21}x_{22}^{2}x_{31}^{2}x_{32}^{2},
x_{11}^{2}x_{12}^{2}x_{21}^{2}x_{22}^{2}x_{31}^{2}x_{32},
x_{11}^{2}x_{12}^{2}x_{21}^{2}x_{22}^{2}x_{31}x_{32}^{2}),
\end{eqnarray*}
and we set $u=x_{11}$.
It follows from Theorem \ref{colon-generalized} that
\begin{eqnarray*}
I_{11}(\mathcal{K}'_{2,2,2}):u
&=&(x_{11}x_{12}x_{21}^{2}x_{22}^{2}x_{31}^{2}x_{32}^{2},
x_{12}^{2}x_{21}^{2}x_{22}^{2}x_{31}^{2}x_{32}^{2},
x_{11}x_{12}^{2}x_{21}^{2}x_{22}x_{31}^{2}x_{32}^{2},\\
&&
x_{11}x_{12}^{2}x_{21}x_{22}^{2}x_{31}^{2}x_{32}^{2},
x_{11}x_{12}^{2}x_{21}^{2}x_{22}^{2}x_{31}^{2}x_{32},
x_{11}x_{12}^{2}x_{21}^{2}x_{22}^{2}x_{31}x_{32}^{2})
\end{eqnarray*}
is again a generalized bi-polymatroidal ideal.
}
\end{Example}

We denote the set of monomial prime
ideals of $T$ by $\mathbf{P}(T)$.
The monomial localization of a monomial ideal $\mathcal{L}\subset T$ with respect
to $\wp\in \mathbf{P}(T)$ is the monomial ideal $\mathcal{L}(\wp)$ obtained from $\mathcal{L}$ by substituting 1 for the variables
that do not belong to $\wp$ \cite{BH}.

Now we want to study the monomial localizations of generalized bi-polymatroidal ideals.
First, we recall the notions that come from the general theory for monomial ideals of $T$ (see for instance \cite{VI}).

\begin{Definition}
\label{log}
Let $T=K[x_{11},\ldots,x_{1m_1},x_{21},\ldots,x_{2m_2},\ldots,x_{n1},\ldots,x_{nm_n}]$ be the polynomial ring  over a field $K$ in the variables
$x_{11},\ldots,x_{1m_1},\ldots,x_{n1},\ldots,x_{nm_n}$.
If $\underline{\Xb}^{(a_1;\dots;a_n)}=x_{11}^{a_{11}} \cdots x_{1m_1}^{a_{1m_{1}}} \cdots x_{n1}^{a_{n1}} \cdots x_{nm_n}^{a_{nm_{n}}}$,
we set
\[
\log(\underline{\Xb}^{(a_1;\dots;a_n)})=(a_1;\dots;a_n)=(a_{11},\dots,a_{1m_{1}},\dots,a_{n1},\dots,a_{nm_{n}})\in \ZZ^{m_{1}+\dots+m_{n}}_{+}.
\]
Given a set $F$ of monomials, the {\em log set} of $F$, denoted by $\log(F)$, consists of all $\log(\underline{\Xb}^{(a_1;\dots;a_n)})$, with
$\underline{\Xb}^{(a_1;\dots;a_n)}\in F$,
\[
\log(F)=\{\log(\underline{\Xb}^{(a_1;\dots;a_n)})=(a_1;\dots;a_n)\in \ZZ^{m_{1}+\dots+m_{n}}_{+}\mid \underline{\Xb}^{(a_1;\dots;a_n)}\in F\}.
\]
\end{Definition}

\begin{Example}
\label{log-example}
{
Let $F=\{x_{11}x_{21}^{2}x_{31}^{2},x_{11}^{3}x_{12}x_{32}^{2},x_{11}x_{31}x_{32}\}$ be a set of monomials in $T=K[x_{11},x_{12},x_{21},x_{22},x_{31},x_{32}]$. The log set of $F$ is
\[
\log(F)=(1,0,2,0,2,0),(3,1,0,0,0,2),(1,0,0,0,1,1)\}.
\]
}
\end{Example}

\begin{Definition}
\label{log-ideal}
Let $\mathcal{L}$ be an ideal of $T$ generated by the set of monomials of $F$. We define
\[
\log(\mathcal{L})=\{\log(\underline{\Xb}^{(a_1;\dots;a_n)})= (a_1;\dots;a_n)\in \ZZ^{m_{1}+\dots+m_{n}}_{+}\mid \underline{\Xb}^{(a_1;\dots;a_n)}\in \mathcal{L}\}.
\]
\end{Definition}

The fundamental result of this section is the following

\begin{Theorem}
\label{localization}
Let $I_{t}(\mathcal{K}'_{m_1,\dots,m_n})$ be the generalized graph ideal of a strong quasi-n-partite graph $\mathcal{K}'_{m_1,\dots,m_n}$.
Then  $I_{t}(\mathcal{K}'_{m_1,\dots,m_n})(\wp)$ is a generalized bi-polymatroidal ideal for all $\wp\in \mathbf{P}(T)$.
\end{Theorem}

\begin{proof}
Suppose that $t, q_{1}, \dots ,q_{n}$ be non negative integers and $\sum_{i=1}^{n}q_{i}=t$, $q_{1}, \dots ,q_{n}\geq 0$. In addition, for each $i=1,\dots,n$ we set $q_{i}\neq t$. Hence
\[
B_{t,2}=\{(a_{1};\dots;a_{n})\in \ZZ^{m_{1}+\dots+m_{n}}_{+}:|a_{i}|=\sum _{j=1}^{m_{i}}a_{ij}=q_{i}, \quad 0\leq a_{ij}\leq 2\}
\]
is the set of bases of a generalized discrete bi-polymatroid of rank $t$.

Let $\underline{\Xb}^{(a_1;\dots;a_n)}\in I_{t}(\mathcal{K}'_{m_1,\dots,m_n})$,
where
$\underline{\Xb}^{(a_1;\dots;a_n)}=\prod_{i=1}^{n}x_{i1}^{a_{i1}} \cdots x_{im_{i}}^{a_{im_{i}}}$,
and
\[
(a_1;\ldots;a_n)=(a_{11},\ldots,a_{1m_1};\ldots;a_{n1},\ldots,a_{nm_n})\in \mathbb{Z}^{m_1+\dots+m_{n}}_{+}.
\]
Then
\[
I_{t}(\mathcal{K}'_{m_1,\dots,m_n})=(\{\underline{\Xb}^{(a_1;\dots;a_n)}\mid  (a_1;\dots;a_n) \in B_{t,2}\}).
\]
For any $ij$, let $(I_{t}(\mathcal{K}'_{m_1,\dots,m_n}))_{\{ij\}}$ be the monomial ideal obtained from
$I_{t}(\mathcal{K}'_{m_1,\dots,m_n})$ by substituting 1 for the variable $x_{ij}$.
This implies that
\[
(I_{t}(\mathcal{K}'_{m_1,\dots,m_n}))_{\{ij\}}=(\{\underline{\Xb}^{(a_1;\dots;a_n)}/x_{ij}^{a_{ij}}|\quad (a_1;\dots;a_n) \in B_{t,2}\}).
\]
We claim that $(I_{t}(\mathcal{K}'_{m_1,\dots,m_n}))_{\{ij\}}$ is a generalized bi-polymatroidal ideal.

The first step is to prove that $(I_{t}(\mathcal{K}'_{m_1,\dots,m_n}))_{\{ij\}}$ is generated in one degree.
If
\[
\mathfrak{r}_{ij}=\max \{a_{ij} \mid (a_1;\dots;a_n)\in B_{t,2}\},
\]
then we prove that
\[
G((I_{t}(\mathcal{K}'_{m_1,\dots,m_n}))_{\{ij\}})=\{\underline{\Xb}^{(a_1;\dots;a_n)}/ x_{ij}^{\mathfrak{r}_{ij}} \mid (a_1;\dots;a_n) \in B_{t,2},\quad a_{ij}=\mathfrak{r}_{ij}\}.
\]
Let $(b_1;\dots;b_n)\in B_{t,2}$; then $b_{ij}< \mathfrak{r}_{ij}$. Now we prove that there exists $(u_1;\dots;u_n) \in B_{t,2}$ with
$u_{ij}=\mathfrak{r}_{ij}$ and such that $\underline{\Xb}^{(u_1;\dots;u_n)}/x_{ij}^{u_{ij}}$ divides $\underline{\Xb}^{(b_1;\dots;b_n)}/x_{ij}^{b_{ij}}$.
We work with induction on $\mathfrak{r}_{ij}-b_{ij}$.

If $\mathfrak{r}_{ij}-b_{ij}=0$, then there is nothing to show.
Suppose that $\mathfrak{r}_{ij}-b_{ij}>0$, i.e., $b_{ij}<\mathfrak{r}_{ij}$.
Let $(a_1;\dots;a_n)\in B_{t,2}$ with $a_{ij}=\mathfrak{r}_{ij}$.
It follows from Theorem \ref{generalized bi-excheng property} that there exists $r$ with $a_{ir}< b_{ir}$ such that
$(a_{1};\ldots;a_{n})-\varepsilon_{ij}+\varepsilon_{ir} \in B_{t,2}.$
Furthermore we set
\[(\alpha_{1};\ldots;\alpha_{n})= (b_{1};\ldots; b_{n})-\varepsilon_{ir}+\varepsilon_{ij}.
\]
Consequently, $\underline{\xb}^{(\alpha_1;\dots;\alpha_n)}/x_{ij}^{\alpha_{ij}}$ divides $\underline{\xb}^{(b_1;\dots;b_n)}/x_{ij}^{b_{ij}}$.
Since
$\mathfrak{r}_{ij}-\alpha_{ij}< \mathfrak{r}_{ij}-b_{ij}$,
the inductive hypothesis shows that there exists $(u_1;\dots;u_n)\in B_{t,2}$ with $u_{ij}=\mathfrak{r}_{ij}$
such that $\underline{\Xb}^{(u_1;\dots;u_n)}/x_{ij}^{u_{ij}}$ divides $\underline{\Xb}^{(\alpha_1;\dots;\alpha_n)}/x_{ij}^{\alpha_{ij}}$.
Therefore, $\underline{\Xb}^{(u_1;\dots;u_n)}/x_{ij}^{u_{ij}}$ divides $\underline{\Xb}^{(b_1;\dots;b_n)}/x_{ij}^{b_{ij}}$, as desired.

The second step is to prove that the set
\[
\mathfrak{B}=\log(G(I_{t}(\mathcal{K}'_{m_1,\dots,m_n})_{\{ij\}}))
\]
is the set of bases of a generalized discrete  bi-polymatroid $\mathcal{Q}$ of rank $t-\mathfrak{r}_{ij}$.
Then $|(\beta_1;\dots;\beta_n)|=t-\mathfrak{r}_{ij}$ for all $(\beta_1;\dots;\beta_n)\in \mathfrak{B}$.

If $(\beta_1;\dots;\beta_n), (\gamma_1;\dots;\gamma_n) \in \mathfrak{B}$ with $\beta_{ir'} > \gamma _{ir'}$, then $r'\neq j$.
By hypothesis for $(a_{1};\ldots; a_{n}), (b_{1};\ldots;b_{n})\in B_{t,2}$, if $a_{ir'}=\beta_{ir'} >\gamma_{ir'}=b_{ir'}$, thus
there exists $r''$ such that $a_{ir''}<b_{ir''}$ and $(a_1;\dots;a_n)-\varepsilon_{ir'}+\varepsilon_{ir''}\in B_{t,2}$.

We set $(\vartheta_{1};\ldots; \vartheta_{n})= (a_{1};\ldots; a_{n})-\varepsilon_{ir'}+\varepsilon_{ir''}$;
since $a_{ij}=b_{ij}=\mathfrak{r}_{ij}$, it follows that $j\neq r''$ and $\vartheta_{ij}=a_{ij}$.
Then $(\vartheta_{1};\ldots; \vartheta_{n})\in \mathfrak{B}$, i.e., $(\beta_{1};\ldots; \beta_{n})-\varepsilon_{ir'}+\varepsilon_{ir''}\in \mathfrak{B}$.
Therefore $\mathfrak{B}$ is a generalized discrete bi-polymatroid, and hence $(I_{t}(\mathcal{K}'_{m_1,\dots,m_n}))_{\{ij\}}$ is a generalized bi-polymatroidal ideal. It follows that $I_{t}(\mathcal{K}'_{m_1,\dots,m_n})(\wp)$ is a generalized bi-polymatroidal ideal for all $\wp\in \mathbf{P}(T)$.
This yields the desired conclusion.
\end{proof}

Corollary \ref{power-bi-polymatroid} with Theorem \ref{localization} now yields

\begin{Corollary}
\label{power-localization}
Let $T=K[x_{11},\ldots,x_{1m_{1}},\ldots,x_{n1},\ldots,x_{nm_n}]$, and let $I_{t}(\mathcal{K}'_{m_1,\dots,m_n})\subset T$ be the generalized graph ideal of a strong quasi-n-partite graph $\mathcal{K}'_{m_1,\dots,m_n}$. Furthermore, let $k$ be a positive integer. Then $I_{t}(\mathcal{K}'_{m_1,\dots,m_n})^{k}(\wp)$ is a generalized bi-polymatroidal ideal for all $\wp\in \mathbf{P}(T)$.
\end{Corollary}

\begin{Example}
\label{L11}
{
Let $T=K[x_{11},x_{12},x_{13},x_{21},x_{22},x_{23}]$ be a polynomial ring and $\mathcal{K}'_{3,3}$ be the strong quasi bi-partite graph on vertices $x_{11},x_{12},x_{13},x_{21},x_{22},x_{23}$.
Let $I_{11}(\mathcal{K}'_{3,3})$ be a generalized graph ideal of $\mathcal{K}'_{3,3}$. Then
\[
I_{11}(\mathcal{K}'_{3,3})= L_{1,5,2}L_{2,6,2}+L_{1,6,2}L_{2,5,2}.
\]
This implies that
\begin{eqnarray*}
&I_{11}(\mathcal{K}'_{3,3})=&
(x_{11}^{2}x_{12}^{2}x_{13}x_{21}^{2}x_{22}^{2}x_{23}^{2},
x_{11}^{2}x_{12}x_{13}^{2}x_{21}^{2}x_{22}^{2}x_{23}^{2},
x_{11}x_{12}^{2}x_{13}^{2}x_{21}^{2}x_{22}^{2}x_{23}^{2},
\\
&&
x_{11}^{2}x_{12}^{2}x_{13}^{2}x_{21}^{2}x_{22}^{2}x_{23},
x_{11}^{2}x_{12}^{2}x_{13}^{2}x_{21}x_{22}^{2}x_{23}^{2},
x_{11}^{2}x_{12}^{2}x_{13}^{2}x_{21}^{2}x_{22}x_{23}^{2}).
\end{eqnarray*}
We put $\wp=(x_{11},x_{12},x_{21},x_{22})$. It then follows that
\begin{eqnarray*}
&I_{11}(\mathcal{K}'_{3,3})(\wp)=&
(x_{11}^{2}x_{12}x_{21}^{2}x_{22}^{2},
x_{11}x_{12}^{2}x_{21}^{2}x_{22}^{2},
x_{11}^{2}x_{12}^{2}x_{21}x_{22}^{2},
\\
&&
x_{11}^{2}x_{12}^{2}x_{21}^{2}x_{22}).
\end{eqnarray*}
Therefore, Corollary \ref{power-localization} yields $I_{11}(\mathcal{K}'_{3,3})^{k}(\wp) $
is a generalized bi-polymatroidal ideal for all $k\geq 1$.
}
\end{Example}

\medskip

\textbf{Data Availability} Data sharing not applicable to this article as no datasets were generated or analysed during the current study.

\medskip

\textbf{Code Availability} Code sharing not applicable to this article as no code was generated the current study.

\medskip

\textbf{Compliance with Ethical Standards}

\medskip

\textbf{Conflict of interests} None


\begin{thebibliography}{plain}

\bibitem{BH}
Bandari, S., Herzog, J.: Monomial localization and polymatroidal ideals, Eur. J. Comb. \textbf{34}, 752-763 (2013)

\bibitem{BM}
Bondy, J. A., Murty, U. S. R.: Graph Theory, GTM \textbf{244}, Springer, (2008)

\bibitem{CH}
Conca, A., Herzog, J.: Castelnuovo-Mumford regularity of products of ideals, Collect. Math. \textbf{54}, 137-152 (2003)

\bibitem{H}
Harary, F.: Graph Theory, Addision Wesley, Reading, MA, (1969)

\bibitem{HHC}
Herzog, J., Hibi, T.: Cohen-Macaulay polymatroidal ideals, European Journal of Combinatorics, \textbf{27}, 513-517 (2006)

\bibitem{HH}
Herzog, J., Hibi, T.: Discrete Polymatroids,  Journal of Alg. Comb. \textbf{16}, 239-268 (2002)

\bibitem{HT1}
Herzog, J., Takayama, Y.:  Resolutions by mapping cones, Homology Homotopy Appl. \textbf{4}, 277-294 (2002)

\bibitem{IB}
Imbesi, M., La Barbiera, M.: Algebraic properties of non-squarefree graph ideals, Math. Rep. \textbf{15}, 107-113 (2013)

\bibitem{BPR}
La Barbiera, M.: A note on unmixed ideals of Veronese bi-type ideal. Turkish Journal of Mathematics, \textbf{37}, 1-7 (2013)

\bibitem{L}
La Barbiera, M.: On standard invariants of bi-polymatroidal ideals.  Algebra Colloq. \textbf{21}, 267-274 (2014)

\bibitem{LM}
La Barbiera, M., Moghimipor, M.: Algebraic and combinatorial properties of generalized bi-polymatroidal ideals, preprint.

\bibitem{VI}
Villareal, H.: Monomial algebras, Second Edition, Monographs and Research Notes in Mathematics, Chapman and Hall/CRC, (2015)


\end{thebibliography}
\end{document}